\newcommand{\pa}{\alpha}
\newcommand{\pb}{\beta}
\newcommand{\pg}{\gamma}
\newcommand{\pd}{\delta}
\newcommand{\pn}{\nu}
\def\pm{\mu}
\newcommand{\pl}{\lambda}
\newcommand{\pk}{\kappa}

\newcommand{\px}{\xi}
\newcommand{\pe}{\varepsilon}
\newcommand{\pf}{\varphi}
\newcommand{\pq}{\theta}
\newcommand{\py}{\psi}

\newcommand{\ph}{\eta}

\newcommand{\po}{\omega}
\newcommand{\proof}{\noindent {\bf Proof.}\hspace{2mm}}

\newcommand{\al}{\aleph}
\newcommand{\rest}{|}
\newcommand{\force}{||\!\!-}

\newcommand{\beg}[1]{\begin{#1}}
\newtheorem{thm}{Theorem}

\newtheorem{lem}[thm]{Lemma}

\newcommand{\se}{\subseteq}
\newcommand{\sm}{\setminus}
\newcommand{\fin}{\raisebox{.6ex}{\framebox[2mm]{}}\par\medskip}

\newcommand{\B}{{\cal B}} 
\renewcommand{\P}{{\cal P}}
\newcommand{\C}{{\cal C}}

\newcommand{\la}{\langle}
\newcommand{\ra}{\rangle}

\newcommand{\cf}{\mbox{\rm cf}}
\newcommand{\dom}{\mbox{\rm Dom}}
\newcommand{\rng}{\mbox{\rm Rng}}
\newcommand{\GMA}{\mbox{\rm GMA}} 
\newcommand{\bc}{\mbox{$\bf c$}}
\newcommand{\N}{{\cal N}_1}
\newcommand{\siz}[1]{\mbox{\bf{$\Sigma^0_{#1}$}}}
\newcommand{\piz}[1]{\mbox{\bf{$\Pi^0_{#1}$}}}
\newcommand{\sio}{\mbox{\bf{$\Sigma^1_1$}}}
\newcommand{\pio}{\mbox{\bf{$\Pi^1_1$}}}
\newcommand{\cub}{\mbox{CUB}}
\newcommand{\nstat}{\mbox{NON-STAT}}
\newcommand{\vx}{\vec{x}}
\newcommand{\vy}{\vec{y}}
\newcommand{\vu}{\vec{u}}
\newcommand{\vv}{\vec{v}}
\newcommand{\vz}{\vec{z}}
\documentclass[12pt]{article}
\author{
Saharon Shelah
\thanks{Research
partially supported by ?
Publication number ?.}\\
Institute of Mathematics\\
Hebrew University\\
Jerusalem, Israel\\
\and
Jouko V\"a\"an\"anen
\thanks{Research partially supported by
grant 1011049 of the Academy of Finland}\\
Department of Mathematics\\
University of Helsinki\\
Helsinki, Finland}
\title{Stationary Sets and Infinitary Logic}
\begin{document}
\maketitle

\begin{abstract}

Let \(K^0_\pl\) be the class of structures
\(\la\pl,<,A\ra\), where \(A\se\pl\) is disjoint from a club,
and let \(K^1_\pl\) be the class of structures
\(\la\pl,<,A\ra\), where \(A\se\pl\) contains a club.
We prove that if \(\pl=\pl^{<\pk}\) is regular,
then no sentence of \(L_{\pl^+\pk}\)
separates \(K^0_\pl\) and \(K^1_\pl\). On the other hand,
we prove that if \(\pl=\pm^+\),
\(\pm=\pm^{<\pm}\), and
a forcing axiom holds (and \(\al_1^L=\al_1\) if \(\pm=\al_0\)),
 then there is a sentence of 
\(L_{\pl\pl}\) which separates \(K^0_\pl\) and \(K^1_\pl\).

\end{abstract}
 


One of the fundamental properties of \(L_{\po_1\po}\) is that
although every countable ordinal itself is definable in  \(L_{\po_1\po}\),
the class of all countable well-ordered structures is not.
In particular, the classes
\beg{eqnarray*}
K^0&=&\{\la\po,R\ra:R\mbox{ well-orders }\po\}\\
K^1&=&\{\la\pl,R\ra:\la\po,R\ra\mbox{ contains a copy of the
rationals}\}
\end{eqnarray*}
cannot be separated
by any \(L_{\po_1\po}\)-sentence. In this paper we consider
infinite quantifier languages \(L_{\pk\pl}\), \(\pl>\po\).
Here well-foundedness is readily definable, but we may instead
consider the class
\[
T_\pl=\{\la\pl,R\ra: \la\pl,R\ra\mbox{ is a tree
with no branches of length \(\pl\)}\}.
\]
If \(\pl=\pl^{<\pl}\), then a result of Hyttinen \cite{hy} 
implies that 
\(T_\pl\) cannot be defined in \(L_{\pl^+\pl}\). 

The main topic of this paper is the question whether the classes 
\beg{eqnarray*}
K^0_\pl&=&\{\la\pl,<,A\ra: A \mbox{ is disjoint from a club of }\pl\}\\
K^1_\pl&=&\{\la\pl,<,A\ra: A \mbox{ contains a club of }\pl\}
\end{eqnarray*}
can be separated in \(L_{\pl^+\pl}\) and related languages.
Note that a set \(A\se\pl\) contains a club if and only if the
tree \(T(A)\) of continuously ascending sequences of elements
of \(A\) has a branch of length \(\pl\).
We show (Theorem~\ref{1st}) that the classes
\(K^0_\pl\) and \(K^1_{\pl}\) cannot
be separated by a sentence of \(L_{\pl^+\pk}\), if
\(\pl=\pl^{<\pk}\) is regular.
The proof of this result uses forcing in a way which seems to be new in the
model theory of infinitary languages.
It follows from this result that  the class
\[S_\pl=\{\la\pl,<,A\ra: A \mbox{ is stationary on }\pl\},\]
that separates \(K^0_\pl\) and \(K^1_{\pl}\), is undefinable
in \(L_{\pl^+\pk}\), if
\(\pl=\pl^{<\pk}\) is regular. 
We complement this result by showing (Theorem~\ref{2nd}) that
if either \(\pl=\pm^+\) and
\(\pm=\pm^{<\pm}>\po\) or \(\pl=\po_1\) and additionally
a forcing axiom holds, then there is a sentence of 
\(L_{\pl\pl}\) which defines \(S_\pl\) and thereby
separates \(K^0_\pl\) and \(K^1_\pl\).

Hyttinen \cite{hy} actually proves more than undefinability of
\(T_\pl\) in \(L_{\pl^+\pl}\). He shows that \(T_\pl\)
is undefinable - assuming \(\pl=\pl^{<\pl}\) - in
\(PC(L_{\pl^+\pl})\). We show (Theorems~\ref{is_pc}
and \ref{no_pc}) that the related statement
that \(S_{\po_1}\) is definable in \(PC(L_{\po_2\po_1})\)
is independent of ZFC+CH.

\section{The case \(\pl=\pl^{<\pm}\).}

\beg{thm}\label{1st} 
If \(\pl=\pl^{<\pk}\) is regular, then 
the classes \(K^0_\pl\) and \(K^1_{\pl}\) cannot
be separated by a sentence of \(L_{\pl^+\pk}\).
\end{thm}

\proof Assume \(\pl=\pl^{<\pk}\) is regular and \(\py\in L_{\pl^+\pk}\).
Let \(\P\) be the forcing notion for adding a Cohen subset to
\(\pl\). Thus \(p\in \P\) if \(p\) is a mapping
\(p:\pa_p\rightarrow 2\) for some \(\pa_p<\pl\).
A condition \(p\) extends another condition \(q\), in symbols
\(p\ge q\), if \(\pa_p\ge\pa_q\) and \(p\rest\-a_q=q\).
Let \(G\) be \(\P\)-generic and \(g=\bigcup G\). Thus
\[V[G]\models g^{-1}(1)\mbox{ is bi-stationary on \(\pl\).}\]
Now either \(\py\) or \(\neg\py\) is true in
\(\langle \pl,<,g^{-1}(1)\rangle\) in \(V[G]\).
We may assume, by symmetry, that it is \(\py\).
Let \(p\in G\) such that 
\[p\force_\P \langle \pl,<,\tilde{g}^{-1}(1)\rangle\models {\py},\]
where \(\tilde{g}\) is the canonical name for \(g\). It is easy
to use \(\pl=\pl^{<\pk}\) and regularity of \(\pl\) to construct
an elementary chain \(\langle M_\px:\px<\pl\rangle\) such that
\begin{itemize}
\item[(i)] \(M_\px\prec\langle H(beth_7(\pl)),\in,<^*\rangle\),
where \(<^*\) is a well-ordering of \(H(beth_7(\pl))\).
\item[(ii)] \(\pl+1\cup\{p\}\cup\{\P\}\cup TC(\{\py\})\subseteq M_0\).
\item[(iii)] \(\langle M_\ph:\ph<\px\rangle\in M_{\px+1}\).
\item[(iv)] \(M_\pn=\bigcup_{\px<\pn}M_\px\) for limit \(\pn\).
\item[(v)] \((M_\px)^{<\pk}\subseteq M_{\px+1}\).
\item[(vi)] \(|M_\px|=\pl\).
\end{itemize}

Let \(M=\bigcup_{\px<\pl}M_\px\). Note, that 
\(M^{<\pk}\subseteq M\) because \(\pl\) is regular.
We shall construct two \(\P\)-generic
sets, \(G^0\) and \(G^1\), over \(M\). For this end, list open
dense \(D\subseteq \P\) with \(D\in M\) as
\(\langle D_\px:\px<\pl\rangle\). Define 
\(G^l=\{p^l_\px:\px<\pl\}\) so that 
\(p^l_0=p\), \(p^l_{\px+1}\ge p^l_\px\) with
\(p^l_{\px+1}\in D_{\px}\cap M\),
\(p^l_{\px+1}(\pa_{p^l_\px})=l\), and
\(p^l_\pn=\bigcup_{\px<\pn}p^l_\px\) for limit
\(\pn\). Clearly, \(G^l\) is \(\P\)-generic over \(M\)
and
\[M[G^l]\models[\langle\pl,<,(g^l)^{-1}(1)\rangle\models\py],\]
where \(g^l=\bigcup G^l\). Note also that 
\(M[G^l]^{<\pk}\subseteq M[G^l]\), because \(M^{<\pk}\subseteq M\) and
\(\P\) is \(<\pk\)-closed.

\begin{lem} If \(\pf(\vx)\in L_{\pl^+\pk}\) such that 
\(TC(\{\pf(\vx)\})\subseteq M\), \(X\in M\), and \(\vec{a}\in \pl^{<\pk}\),
then 
\[\langle\pl,<,X\rangle\models\pf(\vec{a})\iff
M[G^l]\models[\langle\pl,<,X\rangle\models\pf(\vec{a})].\]
\end{lem}

\proof Easy induction on \(\pf(\vx)\). \fin

By the lemma, \(\langle\pl,<,(g^l)^{-1}(1)\rangle\models\py\).
By construction, \(\langle\pl,<,(g^l)^{-1}(1)\rangle\in K^l_\pl\).
Now we can finish the proof. 
Suppose \(K^0_\pl\subseteq Mod(\py)\)
and \(K^1_\pl\cap Mod(\py)=\emptyset\). This contradicts the
fact that \(\langle\pl,<,(g^1)^{-1}(1)\rangle\in K^1_\pl\cap Mod(\py)\).
Suppose \(K^1_\pl\subseteq Mod(\py)\)
and \(K^0_\pl\cap Mod(\py)=\emptyset\). This contradicts 
\(\langle\pl,<,(g^0)^{-1}(1)\rangle\in K^0_\pl\cap Mod(\py)\).
\fin

\beg{cor}
If \(\pl=\pl^{<\pk}\) is regular, then there is no 
\(\pf\in L_{\pl^+\pk}\) such that for all \(A\se\pl\):
\(\langle\pl,<,A\rangle\models\pf\iff
A\mbox{ is stationary}.\)
\end{cor}

Theorem \ref{1st} gives a new proof of the result, referred to above,
that if \(\pl=\pl^{<\pl}\), then \(T_{\pl}\) is not definable
in \(L_{
\pl\pl}\). Our proof does not give the stronger result
that \(T_{\pl}\) is not definable
in \(PC(L_{\pl\pl})\), and there is a good reason: 
\(S_{\po_1}\) may be \(PC(L_{\po_1\po_1})\)-definable, even if
\(2^{\al_0}=\al_1\). This is the topic of the next section.

\section{An application of Canary trees.}

A tree \(\C\) is a {\it Canary tree} if \(\C\) has cardinality 
\(\le 2^{\po}\), \(\C\) has no uncountable branches, but
if a stationary subset of \(\po_1\) is killed by forcing
which does not add new reals, then this forcing adds
an uncountable branch to \(\C\). By \cite{mv}, this is
equivalent to the statement that 
\begin{itemize}
\item[(\(\star\))] For every co-stationary
\(A\se\po_1\) 
there is a mapping \(f\) with \(\rng(f)\se \C\) such that 
for all increasing closed sequences \(s,s'\) of elements of \(A\), if
\(s\) is an initial segment of \(s'\), then \(f(s)<_{\C}f(s')\).
\end{itemize}

\beg{thm}\label{Canary}
\beg{itemize}
\item[(i)] Con(ZF)\(\rightarrow\) Con(ZFC + CH + there is a Canary tree)
\cite{ms}
\item[(ii)] V=L \(\rightarrow\) there are no Canary trees \cite{tv}.
\end{itemize}
\end{thm}

Thus  the non-existence of Canary trees is
consistent with CH, relative to the consistency of ZF. This result
was first proved in \cite{ms} by the method of forcing.

\beg{thm}\label{is_pc}
Assuming CH and the existence of a Canary tree, there is a 
\(\Phi\in PC(L_{\po_2\po_1})\) such that for all \(A\se\po_1\):
\(\langle\po_1,<,A\rangle\models\Phi\iff
A\mbox{ is stationary}.\)
\end{thm} 

\proof Let \(\C\) be a Canary tree.
It is easy to construct a \(PC(L_{\po_2\po_1})\)-sentence \(\Psi\)
such that the following conditions are equivalent
for all \(A\se\po_1\):
\begin{itemize}
\item[(i)] \(\la\po_1,<,A\ra\models\Psi\) 
\item[(ii)] There is a mapping \(f\) with \(\rng(f)\se \C\) such that 
for all increasing closed sequences \(s,s'\) of elements of \(A\), if
\(s\) is an initial segment of \(s'\), then \(f(s)<_{\C}f(s')\).
\end{itemize}
We allow predicate symbols with \(\po\)-sequences of
variables in the \(PC(L_{\po_2\po_1})\)-sentence \(\Psi\).
Now the claim follows from the property (\(\star\)) of Canary trees. \fin

\beg{thm}\label{no_pc}
Con(ZF) implies Con(ZFC + CH + 
there is no
\(\Phi\in PC(L_{\po_2\po_1})\) such that for all \(A\se\po_1\):
\(\langle\po_1,<,A\rangle\models\Phi\iff
A\mbox{ is stationary}\)).
\end{thm}

\proof We start with a model of GCH and add $\al_2$ Cohen subsets to
$\po_1$. In the extension GCH continues to hold. Suppose there is 
in the extension a
\(\Phi\in PC(L_{\po_2\po_1})\) 
such that for all \(A\se\po_1\):
\[\la\po_1,<,A\ra\models\Phi\iff A\mbox{ is stationary}.\]
Since the
forcing to add $\al_2$ Cohen subsets of $\po_1$ satisfies the
$\al_2$-c.c., $\Phi$ belongs to the extension of the universe by $\al_1$
of the subsets. By first adding all but one of the subsets we can work
in $V[A]$ where $A$ is a Cohen subset of $\po_1$ and $\Phi$ is in $V$.
Note that $A$ is a bi-stationary subset of $\po_1$. 
Let
$\P$ be in \(V\) the forcing for adding a Cohen generic subset of $\po_1$ and
let $\tilde{A}$ be the $\P$-name for $A$. 
Let \(p\) force \(\langle\po_1,<,\tilde{A}\rangle\models\Phi\).
By arguing as in the proof of Theorem~\ref{1st}, we can construct
in \(V\)  a model \(M\) of cardinality \(\al_1\) 
containing \(\P\) such that \(M^{\po}\se M\),
\[M\models[p\force\langle\po_1,<,\tilde{A}\rangle\models\Phi],\]
and, furthermore, we can extend \(p\) to a \(\P\)-generic set 
\(H\se\po_1\) over \(M\)
such that  \(H\) is non-stationary.
Thus \(M[H]\) satisfies
\beg{equation}\label{yes}
\langle\po_1,<,H\rangle\models\Phi.
\end{equation}
Now (\ref{yes}) is true in \(V\), because \(M[H]^\po\se M[H]\).
Since \(\P\) is countably closed, we have (\ref{yes}) in
\(V[A]\),
whence \(H\) is stationary in \(V[A]\), contrary
to the fact that \(H\) is non-stationary
in \(V\). \fin 


\section{An application to the topological space \({}^{\po_1}\po_1\).}

Let \(\N\) denote the generalized Baire space consisting of
all functions \(f:\po_1\rightarrow\po_1\),  with the sets
\[N_s=\{f\in\N:f\rest\dom(s)=s\},\]
where \(s\in {}^{<\po_1}\po_1\),
as basic open sets.
We call open sets \(\siz{1}\) and closed sets \(\piz{1}\).
A set of the form \(\bigcup_{\px<\po_1}A_\px\), where
each \(A_\px\) is in \(\bigcup_{\pb<\pa}\piz{\pb}\), is called \(\siz{\pa}\).
Respectively, a set of the form \(\bigcap_{\px<\po_1}A_\px\), where
each \(A_\px\) is in \(\bigcup_{\pb<\pa}\siz{\pb}\), is called \(\piz{\pa}\).
In \(\N\) it is natural to define
Borel sets as follows:
A subset of \(\N\) is {\it Borel}
if it is \(\siz{\pa}\) or \(\piz{\pa}\) for some \(\pa<\po_2\).
A set \(A\se\N\) is \(\pio\) if there is an open set \(B\se
\N\times\N\) such that \(\forall f(f\in A\iff\forall g((f,g)\in B)\).
A set is \(\sio\) if its complement is \(\pio\). 

Let 
\(\cub\) be the set of characteristic functions of closed
unbounded subsets of \(\po_1\), and
\(\nstat\)  the set of characteristic functions of
non-stationary subsets of \(\po_1\). Clearly,
\(\cub\) and \(\nstat\) are disjoint \(\sio\).
It was proved in \cite{mv} that, assuming CH,
\(\cub\) and \(\nstat\) are \(\pio\) if and only if
there is a Canary tree. Another result on \cite{mv}
says that the sets \(\cub\) and \(\nstat\)
cannot be separated by any \(\piz{3}\) 
or \(\siz{3}\) set.

\beg{thm}
Assuming CH, the sets \(\cub\) and \(\nstat\)
cannot be separated by a Borel set.
\end{thm}

\proof 
Let \(\{s_\pa:\pa<\po_1\}\)
enumerate all \(s\in {}^{<\po_1}\po_1\).
Let $ C=\bigcup_{\pa<\po_2}C_\pa$, where
\begin{eqnarray*}
C_{0}&=&\{0,1\}\times\N\\
C_{\pd}&=&\{2,3\}\times{}^{\po_1}(\bigcup_{\pa<\pd}C_\pa).
\end{eqnarray*}
Now we define a Borel set \(B_c\) for each \(c\in C\) as
follows:
\begin{eqnarray*}
B_{(0,f)}=\bigcup_{\pa<\po_1}N_{s_{f(\pa)}}&\mbox{,}&
B_{(1,f)}=\bigcap_{\pa<\po_1}\N\setminus N_{s_{f(\pa)}},\\
B_{(2,f)}=\bigcup_{\pa<\po_1}B_{f(\pa)}&\mbox{,}&
B_{(3,f)}=\bigcap_{\pa<\po_1}B_{f(\pa)}.
\end{eqnarray*}
Clearly, every Borel subset \(X\) of \(\N\) is of the
form \(B_c\) for some \(c\in C\). Then we call
\(c\) a {\it Borel code} of \(X\).

Assume \(A\) is a Borel set 
which separates \(\cub\) and \(\nstat\).
Let \(c\) be a Borel code of \(A\).
Let
\(\P\) be the forcing notion for adding a Cohen subset to
\(\po_1\). 
Let \(G\) be \(\P\)-generic and \(g=\bigcup G\). Thus
\[V[G]\models g^{-1}(1)\mbox{ is bi-stationary.}\]
Now either \(g^{-1}(1)\in B_c\) or \(g^{-1}(1)\in B_c\)
in \(V[G]\).
We may assume, by symmetry, that \(g^{-1}(1)\in B_c\).
Let \(p\in G\) such that 
\[p\force_\P \tilde{g}^{-1}(1)\in B_c,\]
where \(\tilde{g}\) is the canonical name for \(g\). 
Let \(M\prec\langle H(beth_7(\po_1)),\in,<^*\rangle\),
where \(<^*\) is a well-ordering of \(H(beth_7(\pl))\), 
such that 
\(\po_1+1\cup\{p\}\cup\{\P\}\cup TC(\{c\})\subseteq M\),
  \(M^{<{\po_1}}\subseteq M\) and
 \(|M|=\po_1\).

We shall construct two \(\P\)-generic
sets, \(G^0\) and \(G^1\), over \(M\). For this end, list open
dense \(D\subseteq \P\) with \(D\in M\) as
\(\langle D_\px:\px<\po_1\rangle\). Define 
\(G^l=\{p^l_\px:\px<\po_1\}\) so that 
\(p^l_0=p\), \(p^l_{\px+1}\ge p^l_\px\) with
\(p^l_{\px+1}\in D_{\px}\cap M\),
\(p^l_{\px+1}(\pa_{p^l_\px})=l\), and
\(p^l_\pn=\bigcup_{\px<\pn}p^l_\px\) for limit
\(\pn\). Clearly, \(G^l\) is \(\P\)-generic over \(M\)
and
\[M[G^l]\models (g^l)^{-1}(1)\in B_c,\]
where \(g^l=\bigcup G^l\). Note also that 
\(M[G^l]^{<\po}\subseteq M[G^l]\), because \(M^{<\po}\subseteq M\) and
\(\P\) is \(\po\)-closed.

\begin{lem} If \(c\in C\) such that 
\(TC(\{c\})\subseteq M\), and \(f\in M\), then 
\[f\in B_c\iff
M[G^l]\models[f\in B_c].\]
\end{lem}

\proof Easy induction on \(c\). \fin

By the lemma, \((g^l)^{-1}(1)\in B_c\).
By construction, \((g^0)^{-1}(1)\in \nstat\)
and \((g^1)^{-1}(1)\in \cub\).
Now we can finish the proof. 
Suppose \(\cub\subseteq A\)
and \(\nstat\cap A=\emptyset\). This contradicts the
fact that \((g^0)^{-1}(1)\in\nstat\cap A\).
Suppose \(\nstat\subseteq A\)
and \(\cub\cap A=\emptyset\). This contradicts the
fact that \((g^1)^{-1}(1)\in\cub\cap A\).
\fin
\section{The case \(\pl^{\pm}>\pl\).}

Let \(\pm\) be a cardinal. Sets \(A,B\se\pm\) are called
{\it almost disjoint (on \(\pm\))} if \(\sup(A\cap B)<\pm\).
An {\it almost disjoint \(\pl\)-sequence} of subsets of \(\pm\)
is a sequence \(\B=\la B_\pa : \pa<\pl\ra\) such that 
for all \(\pa\ne\pb\), \(|B_\pa|=\pm\) and the sets \(B_\pa\) and
\(B_\pb\) are almost disjoint. The sequence \(\B\) is said to be
{\it definable on} \(L_\pl\) if there is a sequence \(\la\pd_\pa : \pa<\pl\ra\)
such that \(\lim\sup_{\pa<\pl}\pd_\pa=\pl\) and 
the predicates \(x\in B_y\wedge y<\pd_\pa\) and \(x=\pd_y\wedge x<\pa
\wedge y<\pa\) are definable
on every structure \(\la L_\pa,\in\ra\), where \(\pa<\pl\), that is,
there is a first order formula  \(\pf_0(x,y)\)
of the language of set theory such that
for \(x,y<\pa<\pl\): 
\begin{eqnarray*}
x\in B_y\wedge y<\pd_\pa&\iff&\la L_\pa,\in\ra\models\pf_0(x,y).
\end{eqnarray*}

\beg{lem}\label{holds_in_L}
If \(\al_1^L=\al_1\), then there is an almost disjoint
\(\po_1\)-sequence of subsets of \(\po_1\), which
is definable on \(L_{\po_1}\).
\end{lem}

\proof There is a set \(\{B_i:i<\po^L_1\}\)
of almost disjoint subsets of \(\po\)
in \(L\). Since \(\al_1^L=\al_1\),
this set is really of cardinality \(\al_1\).
Let \(\pq(x,y)\) be a \(\Sigma_1\)-formula of set
theory such that for all \(\pa\) and \(x,y\in L_\pa\),
\(x<_L y\iff L_{\pa}\models\pq(x,y)\), where \(<_L\) is the
canonical well-ordering of \(L\).
The claim follows easily. \fin

\beg{thm}\label{2nd}
Suppose 
\beg{itemize}
\item[(i)]\(\pl=\pm^+\).
\item[(ii)] There is an almost disjoint
\(\pl\)-sequence \(\B=\la B_\pa : \pa<\pl\ra\) of subsets
of \(\pm\) which is definable on \(L_\pl\).
\item[(iii)] For all club subsets \(C\) of \(\pl\) there
is a subset \(X\) of \(\pm\) such that for all \(\pa<\pl\)
we have 
\[\pa\in C\iff\sup(B_\pa\sm C)<\pm.\]
\end{itemize}
Then there is a sentence \(\pf\in L_{\pl\pl}\) so that
for all \(A\se\pl\):\[\langle\pl,<,A\rangle\models\pf\iff
A\mbox{ is stationary}.\]
\end{thm}

\proof Suppose \(\pf_0\) defines the almost disjoint
sequence, as above.
We define a sequence of formulas of \(L_{\pl\pl}\). The variable
vectors \(\vx\) in these formulas are always sequences of
the form \(\la x_i:i<\pm\ra\). Let
\(\Phi\) be the conjunction of a large but
finite number of axioms of \(ZFC+V=L\). If \(\psi(\vz)\)
is a formula of set theory, let \(\psi'(\vz,\vx,\vu,\vv)\)
be the result of replacing every quantifier
\(\forall y\ldots\) in \(\Phi\) by \(\forall y(\bigvee_{i<\pm}
y=x_i\rightarrow\ldots)\),
every quantifier
\(\exists y\ldots\) in \(\Phi\) by \(\exists y(\bigvee_{i<\pm}
y=x_i\wedge\ldots)\),
and \(y\in z\) everywhere in \(\Phi\)
by \(\bigvee_{i<\mu}(y=u_i\wedge z=v_i)\).
The following formulas pick \(\pm\) from \(\la\pl,<\ra\):
\[
\begin{array}{lcl}
\pf_{\approx\pm}(y)&\iff& 
\exists \vx((\bigwedge_{i<j<\pm}x_i<x_j)\wedge\forall z
(z<y\leftrightarrow \bigvee_{i<\pm}z=x_i)),\\
\pf_{\in\pm}(y)&\iff&\forall u(\pf_{\approx\pm}(u)\rightarrow y<u),\\
\py_{\in\pm}(\vy)&\iff&\bigwedge_{i<\pm}
\pf_{\in\pm}(y_i)\pf_{B,1}(x,\vu,\vv,z,y)
\end{array}\]
The following formulas are needed to refer to well-founded
models of set theory:
\[
\begin{array}{lcl}
\pf_{uni}(\vx,z)&\iff&\bigvee_{i<\pm}z=x_i\\
\pf_{eps}(\vx,\vu,\vv,z,y)&\iff&\pf_{uni}(\vx,z)\wedge \pf_{uni}(\vx,y)
\wedge \bigvee_{i<\pm}(z=u_i\wedge y=v_i)\\
\pf_{wf}(\vx,\vu,\vv)&\iff&\Phi'(\vx,\vu,\vv)
\wedge\forall\vy
((\bigwedge_{i<\pm}\pf_{uni}(\vx,y_i))\rightarrow\\
&&\bigvee_{i<\pm}\neg 
\pf_{eps}(\vx,\vu,\vv,y_{i+1},y_i))\\
\pf_{cor}(\vx,\vu,\vv,z)&\iff&\forall s(s<z\leftrightarrow
\bigvee_{i<\pm}(s=u_i\wedge z=v_i)\\
\end{array}\]
Let
\[
\begin{array}{lcl}
\pf_B(z,y)&\iff&\exists \vx\exists\vu\exists\vv
(\pf_{wf}(\vx,\vu,\vv)
\wedge\pf_{cor}(\vx,\vu,\vv,z)\wedge\\
&&\pf_{cor}(\vx,\vu,\vv,y)\wedge\phi_0'(z,y,\vx,\vu,\vv)).\\
\end{array}\]
The point is that if  \(\pa\in\pm\) and \(\pb\in \pl\), then
\(\pa\in B_\pb\) if and only if
\(\la\pl,<\ra\models \pf_B(\pa,\pb)\).
The following formula says that the element \(y\)
of \(\pm\) is in the subset of \(\pl\) coded
by \(\vx\):
\[
\begin{array}{lcl}
\pf_\pe(y,\vx)&\iff&\exists u(\pf_{\in\pm}(u)\wedge\forall z((
\pf_B(z,y)\wedge\bigwedge_{i<\pm}z\ne x_i)\rightarrow z<u),\\
\end{array}\]
Finally, if:
\[\begin{array}{lcl}
\pf_{ub}(\vx)&\iff&\forall y\exists z(y<z\wedge
\pf_\pe(z,\vx)),\\
\pf_{cl}(\vx)&\iff&\forall y
(\forall z(z<y\rightarrow\exists u(z<u\wedge u<y\wedge\pf_\pe(u,\vx)))
\rightarrow\\
&&\pf_\pe(y,\vx))\\
\pf_{cub}(\vx)&\iff&\pf_{ub}(\vx)\wedge\pf_{cl}(\vx)\\
\pf_{stat}&\iff&\forall\vx((\py_{\in\pm}(\vx)\wedge
\pf_{cub}(\vx))\rightarrow\exists y(A(y)\wedge\pf_\pe(y,\vx))),\\
\end{array}\] 
then \(\la\pl,<,A\ra\models\pf_{stat}\)  if and only if 
\(A\) is stationary.
\fin

\beg{cor}\label{special_positive}
If \(2^{\al_0}>\al_1\), \(\al_1^L=\al_1\)
and MA, then there is a \(\pf\in L_{\po_1\po_1}\)
such
that for all \(A\se\po_1\):\[\langle\po_1,<,A\rangle\models\pf\iff
A\mbox{ is stationary}.\]
\end{cor}

\proof We choose \(\pl=\po_1\) and \(\pm=\po_0\) in 
Theorem~\ref{2nd}. Condition (ii) holds by Lemma~\ref{holds_in_L}.
Condition (iii) is a consequence of MA + $\neg$CH
by \cite{msol}. \fin

\medskip

{\bf Note.} The proof of Corollary~\ref{special_positive} shows that
we actually get the following stronger result: 
If \(2^{\al_0}>\al_1\), \(\al_1^L=\al_1\)
and MA, then the full second order extension \(L^{II}_{\po_1\po_1}\)
of \(L_{\po_1\po_1}\) is reducible to \(L_{\po_1\po_1}\) in
expansions of \(\la\po_1,<\ra\). Then, in particular,
\(T_{\al_1}\) is \(PC(L_{\po_1\po_
1})\)-definable. 
This kind of reduction cannot hold
on all models. For example, \(\po_1\)-like dense linear orders
with a first element are all \(L_{\infty\po_1}\)-equivalent, but not 
\(L^{II}_{\po\po}\)-equivalent.

\medskip

For \(\pa<\pl=\pm^+\), let \(\la a^\pa_i:i<\pm\ra\)
be a continuously increasing sequence of subsets of \(\pa\) with
\(\pa=\bigcup_{i<\pm}a^\pa_i\) and \(|a^\pa_i|<\pm\).
Define \(f_\pa:\pm\rightarrow\pm\) by
\[f_\pa(i)=otp(a^\pa_i).\]
Let \(D_\pm\) be the club-filter on \(\pm\).
Define for \(f,g\in {}^\pm\pm\);
\[f\sim_{D_\pm}g\iff\{i:f(i)=g(i)\}\in D_\pm.\]

\beg{lem}\label{independence}
\(f_\pa/D_\pm\) is independent of the
choice of the sequence  \(\la a^\pa_i:i<\pm\ra\).
\end{lem}

\beg{thm}\label{mu_big}
Suppose 
\beg{itemize}
\item[(i)] \(\pl=\pm^+\), where \(\pm=\pm^{<\pm}>\al_0\).
\item[(ii)] For every club \(C\se\pl\) there is some \(X\se\pm\times\pm\)
such that 
\beg{eqnarray*}
\pa\in C&\rightarrow&
\{i<\pm:(i,f_\pa(i))\in X\}\mbox{ contains a club}\\
\pa\not\in C&\rightarrow&
\{i<\pm:(i,f_\pa(i))\not\in X\}\mbox{ contains a club}.
\end{eqnarray*}
\end{itemize}
Then there is a sentence \(\pf\in L_{\pl\pl}\) such that
for all \(A\se\pl\):\[\langle\pl,<,A\rangle\models\pf\iff
A\mbox{ is stationary}.\]
\end{thm}

\proof This is like the proof of Theorem~\ref{2nd}. One uses
Lemma~\ref{independence} to refer to the functions \(f_\pa\). 
We leave the details to the reader. \fin

The {\it Generalized Martin's Axiom for \(\pm\)} (\(\GMA_\pm\)) from
\cite{sh}
is the following principle:
\[\mbox{
\beg{minipage}[t]{12cm}
\noindent Suppose \(\P\) is a forcing notion with the properties:
\beg{description}
\item[(GMA1)] Every descending sequence of length \(<\pm\) in \(\P\)
has a greatest lower bound.
\item[(GMA2)] If \(p_\pa\in\P\) for \(\pa<\pm^+\), then
there is a club \(C\se\pm^+\) and a regressive function
\(f:\pm^+\rightarrow\pm^+\) such that
if \(\pa\in C\) and \(\cf(\pa)=\pm\), then the set
\[A=\{p_\pb:\cf(\pb)=\pm,f(\pa)=f(\pb)\}\]
is well-met (i.e. \(p,q\in A \rightarrow p\vee q\in a\)).
\end{description}
Then for any dense open sets \(D_\pa\se\P\), \(\pa<\pk\),
where \(\pk<2^{\pm}\), there is a filter in \(\P\) which
meets every \(D_\pa\).
\end{minipage}}\]

\beg{pro}\label{GMA}
Suppose \(\pl=\pm^+\), where \(\pm=\pm^{<\pm}>\al_0\), and \(\GMA_\pm\).
Then for every club \(C\se\pl\) there is some \(X\se\pm\times\pm\)
such that 
\beg{eqnarray*}
\pa\in C&\rightarrow&
\{i<\pm:(i,f_\pa(i))\in X\}\mbox{ contains a club}\\
\pa\not\in C&\rightarrow&
\{i<\pm:(i,f_\pa(i))\not\in X\}\mbox{ contains a club}.
\end{eqnarray*}
\end{pro}

\proof Let a club \(C\se\pl\) be given. For \(\pa<\pb<\pl\), let
\(C_{\pa\pb}\in D_\pm\) so that 
\(f_\pa\rest C_{\pa\pb}<f_\pb\rest C_{\pa\pb}\). Let \(\P\) consist
of conditions
\[p=(B^p,f^p,\bc^p,g^p),\]
where
\beg{itemize}
\item[(i)] \(B^p\se\pl\). \(|B^p|<\pm\).
\item[(ii)] \(f^p\) is a partial mapping with
\(\dom(f^p)\se\pm\times\pm\), 
\(|\dom(f^p)|<\pm\), and \(\rng(f^p)\se\{0,1\}\).  
\item[(iii)] If \(\pa\in B^p\), then 
\(\{i<\pm : (i,f_\pa(i))\in\dom(f^p)\}\)
is an ordinal \(j^p_\pa\).
\item[(iv)] \(\bc^p=\la c^p_\pa:\pa\in B^p\ra\), where \(c^p_\pa\) is a closed
subset of \(j^p_\pa\). We denote \(\max(c^p_\pa)\) by \(\pd^p\).
\item[(v)] If \(\pa\in B^p\cap C\) and \(i\in C^p_\pa\), then
\(f^p(i,f_\pa(i))=1\).
If \(\pa\in B^p\sm C\) and \(i\in C^p_\pa\), then
\(f^p(i,f_\pa(i))=0\).
\item[(vi)] \(g^p\) is a partial mapping with 
\(\dom(g^p)\se[B^p]^2\) and \(\rng(g^p)\se\pm\).
\item[(vii)] If \(\pa< \pb\in\dom(g^p)\), then
\(\emptyset\ne c^P_\pa\sm g(\pa,\pb)\se C_{\pa\pb}\).
\end{itemize}
The partial ordering ``\(q\) extends \(p\)'' is defined as
follows:
\beg{eqnarray*}
p\le q &\Leftrightarrow&
B^p\se B^q, f^p\se f^q, g^p\se g^q,\\
&&\forall\pa\in B^p(c^p_\pa\mbox{ is an initial segment of \(c^q_\pa\)}),\\
&&\mbox{ and if \(\pd^p<\pd^q\), then \(\dom(g^q)\supseteq
[B^p]^2\)}.
\end{eqnarray*}
We show now that \(\P\) satisfies conditions (GMA1) and (GMA2).

\beg{lem}\label{GMA1}
\(\P\) satisfies (GMA1).
\end{lem}

\proof
Let \(p_o\le\ldots\le p_i\le\ldots (i<\pg)\) in \(\P\) with \(\pg<\pm\).
We may assume \(\pd^{p_0}<\pd^{p_1}<\ldots\). Let
\(\pd=\sup\{\pd^{p_i}:i<\pg\}\).
Let \(B=\bigcup_{i<\pg}B^{p_i}\). We extend \(\bigcup_if^{p_i}\)
to \(f\) by defining
\[f(\pd,f_\pa(\pd))=\left\{
\begin{array}{ll}
1&\mbox{ if \(\pa\in B\cap C\)}\\
0&\mbox{ if \(\pa\in B\sm C\).}\end{array}\right.\]
We have to check that this definition is coherent, i.e.,
if \(\pa\in B\cap C\) and \(\pb\in B\sm C\), then 
\(f_\pa(\pd)\ne f_\pb(\pd)\). Suppose
\(\pa\in B^{p_i}\) and \(\pb\in B^{p_{i'}}\) with 
\(\pa<\pb\) and \(i<i'\). Since \(\pd^{p_i}<\pd^{p_{i'}}\),
\(g(\pa,\pb)\) is defined and \(c^{p_i}_\pa\sm g(\pa,\pb)\se
C_{\pa\pb}\). Hence \(\pd\in C_{\pa\pb}\),
whence \(f_\pa(\pd)<f_\pb(\pd)\).

Let \(\bc=\la c_\pa:\pa\in B\ra\) 
where \(c_\pa=\bigcup_i c^{p_i}_\pa\cup\{\pd\}\).
Let \(j=\bigcup_i j^{p_i}\cup\{\pd\}\).
Now the condition \(p=(B,f,\bc,g)\) is the
needed l.u.b. of \((p_i)_{i<\pm}\). \fin

\beg{lem}\label{GMA2}
\(\P\) satisfies (GMA2).
\end{lem}

\proof
Suppose \(p_\alpha\), \(\alpha<\pl\), are in \(\P\).
Let \(h\) be a one-one mapping from \(\P\) to 
odd ordinals \(<\pl\). By
\(\pm^{<\pm}=\mu\) there is a club \(C\se\pl\)
such that if \(\pa\in C\), \(\cf(\pa)=\pm\),
and \(B^p\se\pa\), then \(h(p)<\pa\),
and if \(\pa<\pb\), \(\pa,\pb\in C\), then \(B^{p_\pa}\se\pb\). 
Choose a regressive function \(g\) from the complement of 
\(C\) to the even ordinals that is one-one on ordinals of cofinality
\(\pm\). Suppose \(\cf(\pa)=\pm\). Let \(f(\pa)=g(\pa)\)
if \(\pa\not\in C\), and \(f(\pa)=h(p_\pa\rest\pa)\) if
\(\pa\in C\). Suppose now \(\pa<\pb\),
\(\cf(\pa)=\cf(\pb)=\pm\),
and \(f(\pa)=f(\pb)\). W.l.o.g. \(\pa,\pb\in C\).
Thus \(h(p_\pa\rest\pa)=h(p_\pb\rest\pb)\),
whence \(p_\pa\rest\pa = p_\pb\rest\pb\).
It follows that \(p_\pa\) and \(p_\pb\) have a l.u.b.
\fin

Let \[D_{\pa\pb}=\{p\in\P:\pa\in B^p\mbox{ and }\pd^p\ge\pb\}\]
where \(\pa<\pl\), \(\pb<\pm\). We show that 
\(D_{\pa\pb}\) is dense open. Suppose therefore \(p\in\P\)
is given. We construct \(q\in D_{\pa\pb}\) with
\(p\le q\). Let \(B^q=B^p\cup\{\pa\}\).
Let
\[E=\bigcap\{C_{\px\ph}: \px,\ph\in B^q, \px<\ph\} (\in D_{\pm}).\]
Let \(\pd^q\in E\setminus\pb\). Define 
\(\bc^q=\la c^q_\px:\px\in B^q\ra\) by
\[c^q_\px=\left\{\begin{array}{ll}
c^p_\px\cup\la\pd^q\ra,&\mbox{ if }\px\ne\pa\\
\la\pd^q\ra,&\mbox{ if }\px=\pa.
\end{array}\right.\]
Let
\[f^q=f^p\cup\left\{\begin{array}{ll}
\{((j,f_\pa(j)),1):j^p\le j\le\pd^q\},&\mbox{ if }\pa\in C\\
\{((j,f_\pa(j)),0):j^p\le j\le\pd^q\},&\mbox{ if }\pa\not\in C.
\end{array}\right.\]
Let \(g^q(\px,\ph)=\pd^p\) for \((\px,\ph)\in
[M^p]^2\setminus\dom(g^p)\).
Let \(q=(B^q,f^q,g^q,\pd^q)\).
Then \(q\in D_{\pa\pb}\), and \(p\le q\).

Let \(G\) be a filter that meets every \(D_{\pa\pb}\).
Let 
\beg{eqnarray*}
B&=&\bigcup\{B^p:p\in G\}\\
f&=&\bigcup\{f^p:p\in G\}\\
c_\pa&=&\bigcup\{c^p_\pa:p\in G\}
\end{eqnarray*}
Then \(B=\pl\) and each \(c_\pa\) is a club of
\(\pm\). Let \(X=\{(\pa,\pb)\in\pm\times\pm:
f(\pa,\pb)=1\}\). Suppose \(\pa\in C\) and \(i\in c_\pa\).
Then \(f(i,f_{\pa}(i))=1\) whence \((i,f_\pa(i))\in X\).
Suppose \(\pa\not\in C\) and \(i\in c_\pa\).
Then \(f(i,f_\pa(i))=0\) whence
\((i,f_\pa(i))\not\in X\). \fin

\beg{cor}
Suppose \(\pl=\pm^+\), where \(\pm=\pm^{<\pm}>\al_0\), and \(\GMA_\pm\).
Then there is a sentence \(\pf\in L_{\pl\pl}\) such that
for all \(A\se\pl\):\[\langle\pl,<,A\rangle\models\pf\iff
A\mbox{ is stationary}.\]
\end{cor}

\proof
The claim follows from Theorem~\ref{mu_big} and 
Proposition~\ref{GMA}. \fin

\end{document}